\documentclass[12pt]{article}
\usepackage{a4}
\usepackage{amsfonts}
\usepackage{amssymb}

\begin{document}

\title{\textbf{The asymptotic behaviour of the number of solutions of polynomial congruences}}
\author{Dirk Segers\thanks{Postdoctoral Fellow of the
Fund for Scientific Research - Flanders (Belgium). \newline The
author is supported by FWO-Flanders project G.0318.06. \newline
}}

\date{August 23, 2012}

\maketitle

\begin{abstract}
One mentions in a lot of papers that the poles of Igusa's $p$-adic
zeta function determine the asymptotic behavior of the number of
solutions of polynomial congruences. However, no publication
clarifies this connection precisely. We try to get rid of this
gap.
\end{abstract}

\section{Introduction}

\noindent \textbf{(1.1)} Let $f \in \mathbb{Z}[x_1,\ldots,x_n]$ be
a polynomial over the integers in $n$ variables. Put
$x=(x_1,\ldots,x_n)$. We are interested in the number of solutions
of $f(x) \equiv 0 \mbox{ mod } m$ in $(\mathbb{Z}/m\mathbb{Z})^n$
for an arbitrary positive integer $m$. The Chinese remainder
theorem reduces this problem to the case that $m$ is a power of a
prime. Fix from now on a prime $p$ and let $M_i$, $i \in
\mathbb{Z}_{\geq 0}$, be the number of solutions of the congruence
$f(x) \equiv 0 \mbox{ mod } p^i$ in
$(\mathbb{Z}/p^i\mathbb{Z})^n$. The aim of this paper is to study
the asymptotic behaviour of the numbers $M_i$, and to relate this
behaviour with information about the poles of Igusa's $p$-adic
zeta function, which will be defined in (1.3).

\vspace{0,5cm}

\noindent \textbf{(1.2)} Let $\mathbb{Z}_p$ be the ring of
$p$-adic integers. The behaviour of the $M_i$ is well understood
if $f^{-1}\{0\}$ has no singular point in $\mathbb{Z}_p^n$.
Indeed, we can take a $k \in \mathbb{Z}_{>0}$ for which $f$ has no
singular point modulo $p^k$ because $f$ has no singular point in
the sequentially compact space $\mathbb{Z}_p^n$. Using Hensels
lemma, one obtains that $M_i = M_{2k-1} p^{(n-1)(i-2k+1)}$ for
every $i \geq 2k-1$.

\vspace{0,5cm}

\noindent \textbf{(1.3)} The behaviour of the $M_i$ is more
complicated if $f$ has a singular point in $\mathbb{Z}_p^n$. At
this stage, we introduce Igusa's $p$-adic zeta function $Z_f(s)$
of $f$. It is defined by
\[ Z_f(s)= \int_{\mathbb{Z}_p^n} |f(x)|^s \, |dx| \]
for $s \in \mathbb{C}$, $\mbox{Re}(s) > 0$, where $|dx|$ denotes
the Haar measure on $\mathbb{Q}_p^n$, so normalized that
$\mathbb{Z}_p^n$ has measure $1$. Note that $Z_f(s)$ only depends
on $p^{-s}$. We will write $Z_f(t)$ if we consider $Z_f(s)$ as a
function in the variable $t:=p^{-s}$.

All the $M_i$ describe and are described by $Z_f(t)$ through the
equivalent relations
\[
Z_f(t) = P(t) - \frac{P(t)-1}{t} \qquad \mbox{and} \qquad P(t) =
\frac{1-t Z_f(t)}{1-t},
\]
where the Poincar\'e series $P(t)$ of $f$ is defined by
\[
P(t) = \sum_{i=0}^{\infty} M_i(p^{-n}t)^i.
\]

\vspace{0,5cm}

\noindent \textbf{(1.4)} Igusa proved in \cite{Igusahigherdegree}
that $Z_f(s)$ is a rational function of $p^{-s}$ by calculating
the integral on an embedded resolution of the singularities of
$f$, which always exists by Hironaka's theorem \cite{Hironaka}.
This implies that it extends to a meromorphic function $Z_f(s)$ on
$\mathbb{C}$, which is also called Igusa's $p$-adic zeta function
of $f$. We also obtain from the relations in (1.3) that $P(t)$ is
a rational function.

Igusa determined actually a specific form of the rational function
which involves geometric data of an embedded resolution $g$ of $f$.
He obtained that $Z_f(t)$ can be written in the form
\[
Z_f(t) = \frac{A(t)}{\prod_{j \in J} (1-p^{-\nu_j}t^{N_j})},
\]
where $A(t) \in S[t]$, with $S := \{z/p^i \mid z \in \mathbb{Z},i
\in \mathbb{Z}_{\geq 0} \}$, where $A(t)$ is not divisible by any
of the $1-p^{-\nu_j}t^{N_j}$ and where the $N_j$ and $\nu_j-1$ are
the multiplicities of $f \circ g$ and $g^*dx$ along an irreducible
component $E_j$ of $g^{-1}(f^{-1}\{0\})$. It is surprising that
most irreducibele components of $g^{-1}(f^{-1}\{0\})$ do not
induce a factor in the denominator. This would be elucidated if
the monodromy conjecture (see for example \cite{Denefreport}) is
true.

It follows from (1.3) and $Z_f(t=1)=1$ that we can write
\[
P(t) = \frac{B(t)}{\prod_{j \in J} (1-p^{-\nu_j}t^{N_j})},
\]
where $B(t) \in S[t]$. Here, $B(t)$ is also not divisible by any of
the $1-p^{-\nu_j}t^{N_j}$. The poles of $P(t)$ and $Z_f(t)$ are
actually the same.

\vspace{0,5cm}

\noindent \textbf{(1.5)} In this paper, we try to explain the
relation between the poles (and their order) of $P(t)$, which are
the same as those of $Z_f(t)$, and the numbers $M_i$. If also the
principal parts of the Laurent series of $P(t)$ at all poles are
known, we will even calculate the numbers $M_i$ (and not only
their asymptotic behaviour) for $i$ large enough. The principal
parts of the Laurent series of $Z_f(t)$ and $P(t)$ at a certain
pole can be calculated from each other, which is also the case for
the ones of $Z_f(s)$ and $Z_f(t)$ at corresponding poles.
Therefore, it is also possible to calculate the numbers $M_i$ for
$i$ large enough from the principal parts of the Laurent series of
$Z_f(s)$ at all its poles. This will not be worked out in the
paper because it leads to formulas which are more complicated and
which do not give us more insight.

\vspace{0,5cm}

\noindent \textsl{Reference.} An introduction to Igusa's $p$-adic
zeta function which contains more details can be found in
\cite[Section 1.1]{Segersthesis}, \cite{Igusabook} or
\cite{Denefreport}.

\vspace{0,5cm}

\noindent \textsl{Acknowledgements.} I want to thank Pierrette
Cassou-Nogu\`es for pointing my attention at this problem.

\section{The asymptotic behaviour}

\noindent \textbf{(2.1)} We define an equivalence relation on $J$.
We say that $j_1 \sim j_2$ iff $\nu_{j_1}/N_{j_1} =
\nu_{j_2}/N_{j_2}$. This equivalence relation determines a
partition of $J$ into sets $J_k$, $k \in V$. For $k \in V$, we
denote the lowest common multiple of the $\nu_j$, $j \in J_k$, by
$a_k$ and the lowest common multiple of the $N_j$, $j \in J_k$, by
$b_k$. Remark that $a_k/b_k = \nu_j/N_j$ for all $j \in J_k$. Let
$m_k$ be the cardinality of $J_k$. Because $1-p^{-a_k}t^{b_k}$ is
a multiple of $1-p^{-\nu_j}t^{N_j}$ for all $j \in J_k$, we can
write
\[
P(t) = \frac{C(t)}{\prod_{k \in V} (1-p^{-a_k}t^{b_k})^{m_k}},
\]
where $C(t) \in S[t]$.

\vspace{0,5cm}

\noindent \textbf{Theorem.} \textsl{There exists a unique
decomposition of every $M_i$ with $i > \deg(P(t))$ of the form
\[
M_i = \sum_{k=1}^r g_k(i)p^{\ulcorner l_k i \urcorner},
\]
where the $l_k$ are different rational numbers and where every
$g_k(i)$ is a nonzero function which is polynomial with rational
coefficients on residue classes. The $r$ numbers $l_k-n$ are the
real parts of the poles of $Z_f(s)$. If we denote the elements of
$V$ by $1,\ldots,r$ in such a way that $l_k-n = -a_k/b_k$ for
every $k \in \{1,\ldots,r\}$, we have for $k \in \{1,\ldots,r\}$
that
\begin{enumerate}
\item the function $g_k(i)$ is polynomial on each residue class modulo
$b_k$,
\item the maximum of the degrees of these polynomials is equal to
$m_k-1$ and
\item these polynomials (and thus also $g_k(i))$ are determined by
the principal parts of the Laurent series of $Z_f(s)$ in the poles
with real part $-a_k/b_k$.
\end{enumerate}}

\vspace{0,2cm}

\noindent \textsl{Remark.} (1) The $l_k$ are rational numbers less
than $n$ because the real parts of the poles of $Z_f(s)$ are
negative rational numbers. The author proved in
\cite{Segersmathann} that the real part of every pole of $Z_f(s)$
is larger than or equal to $-n/2$. This implies that $l_k \geq
n/2$ for every $k \in \{1,\ldots,r\}$. Moreover, in the case that
$n=3$ and $f$ has no singular point in $\mathbb{Z}_p^3$ of
multiplicity 2, the author proved \cite{Segersmathz} that there
are no poles with real part less than $-1$, which implies that
$l_k \geq 2$ for every $k \in \{1,\ldots,r\}$.

(2) It follows from the theorem that the asymptotic behaviour of
the number of solutions is determined by the largest real part of
a pole of $Z_f(s)$ and by the largest order of a pole with maximal
real part.

\vspace{0,2cm}

\noindent \textsl{Proof.} Applying decomposition into partial
fractions in $\mathbb{Q}[t]$, we can write
\[
P(t) = C_0(t) + \sum_{k \in V}
\frac{C_k(t)}{(1-p^{-a_k}t^{b_k})^{m_k}},
\]
where every $C_k(t) \in \mathbb{Q}[t]$ and where $\deg(C_k(t)) < m_k
b_k$ for $k \in V$. Note that the term $C_0(t)$ does not give a
contribution to $M_i$ for $i > \deg(C_0(t))$ and that $\deg(C_0(t))
= \deg(P(t))$ if one of them is non-negative. Now we look at the
contributions of the other terms. So fix $k \in V$. Note that
$C_k(t)$ contains exactly the information of the principal parts of
the Laurent series of $P(t)$ at the poles with absolute value
$p^{a_k/b_k}$. We have
\begin{eqnarray*}
\frac{C_k(t)}{(1-p^{-a_k}t^{b_k})^{m_k}} & = &
\frac{C_{k,m_k}(t)}{(1-p^{-a_k}t^{b_k})^{m_k}} +
\frac{C_{k,m_k-1}(t)}{(1-p^{-a_k}t^{b_k})^{m_k-1}} + \cdots +
\frac{C_{k,1}(t)}{1-p^{-a_k}t^{b_k}} \\ & = & \sum_{d=0}^{b_k-1}
\sum_{e=0}^{\infty} g_{k,d}(e) p^{-ea_k} t^{eb_k+d} \\ & = &
\sum_{d=0}^{b_k-1} \sum_{e=0}^{\infty} g_{k,d}(e) p^{\llcorner
da_k/b_k \lrcorner} p^{\ulcorner (n-a_k/b_k)(eb_k+d) \urcorner}
\frac{t^{eb_k+d}}{p^{n(eb_k+d)}},
\end{eqnarray*}
where $C_{k,l}(t) \in \mathbb{Q}[t]$ with $\deg(C_{k,l}(t))<b_k$ and
where the maximum of the degrees of the polynomials $g_{k,d}(e)$ is
equal to $m_k-1$. Actually, if we denote the coefficient of $t^d$ in
$C_{k,l}(t)$ by $C_{k,l,d}$, we get
\[
g_{k,d}(e) = C_{k,m_k,d} \frac{(e+m_k-1)!}{e!(m_k-1)!} +
C_{k,m_k-1,d} \frac{(e+m_k-2)!}{e!(m_k-2)!} + \cdots +
C_{k,1,d}.
\]
\begin{flushright}
$\Box$
\end{flushright}

\vspace{0,5cm}

\noindent \textbf{(2.2)} Finally, we give two examples. In the first
example, all the coefficients of the polynomials $C_k(t)$, $k \in
V$, are in $S$. This is in some sense the easiest situation. The
second example shows that this is not always the case. There are
several ways to compute the Poincar\'e series: one can calculate the
integral on an embedded resolution of singularities of $f$, one can
use the formula for polynomials which are non-degenerated over
$\mathbb{F}_p$ with respect to their Newton polyhedron
\cite{DenefHoornaert} and one can use the $p$-adic stationary phase
formula \cite[Theorem 10.2.1]{Igusabook}. All these techniques are
also explained in \cite[Section 1.1]{Segersthesis}.

\vspace{0,5cm}

\noindent \textsl{Example 1.} Let $f(x,y) = y^2-x^3$ and let $p$ be
an arbitrary prime. Then,
\begin{eqnarray*}
P(t) & = & \frac{-t^6+p^4t^2-p^3t^2+p^6}{(p^5-t^6)(p-t)} \\ & = &
\frac{2p^{-5}t^5 + 2p^{-4}t^4 + 2p^{-3}t^3 + 2p^{-2}t^2 +
(p+1)p^{-2}t + (p+1)p^{-1}}{1-p^{-5}t^6} - \frac{p^{-1}}{1-p^{-1}t}.
\end{eqnarray*}
We obtain for every $e \in \mathbb{Z}_{\geq 0}$ that
\[ \begin{array}{rclcccrcl}
M_{6e} & = & (p+1)p^{7e-1} - p^{6e-1}, &  &  &  &  M_{6e+1} & = &
(p+1)p^{7e} - p^{6e}, \\  M_{6e+2} & = & 2p^{7e+2}-p^{6e+1}, &  &
& & M_{6e+3} & = & 2p^{7e+3}-p^{6e+2}, \\ M_{6e+4} & = &
2p^{7e+4}-p^{6e+3} &  &  \mbox{and} &  &   M_{6e+5} & = &
2p^{7e+5}-p^{6e+4}.
\end{array} \]

\vspace{0,5cm}

\noindent \textsl{Example 2.} Let $f(x,y) = x^3 + y^5$ and let $p$
be an arbitrary prime. Then,
\begin{eqnarray*}
P(t) & = & \frac{\begin{array}{c} -t^{15} + (p-1)t^{14} +
(p-1)pt^{12} + (p-1)p^3t^9 \\ + (p-1)p^3t^8 + (p-1)p^5t^5 +
(p-1)p^6t^3 + (p-1)p^6t^2 + p^9 \end{array}}{(p^8-t^{15})(p-t)} \\
& = & \frac{C_1(t)}{1-p^{-8}t^{15}} + \frac{C_2(t)}{1-p^{-1}t},
\end{eqnarray*}
where
\begin{eqnarray*}
C_1(t) & = & \frac{p^7+p-2}{(p^7-1)p^8}t^{14} +
\frac{p^7+p^2-p-1}{(p^7-1)p^8}t^{13} +
\frac{p^7+p^2-p-1}{(p^7-1)p^7}t^{12} +
\frac{p^7+p^3-p^2-1}{(p^7-1)p^7}t^{11} \\ & & +
\frac{p^7+p^3-p^2-1}{(p^7-1)p^6}t^{10} +
\frac{p^7+p^3-p^2-1}{(p^7-1)p^5}t^9 +
\frac{p^7+p^4-p^3-1}{(p^7-1)p^5}t^8 +
\frac{p^7+p^5-p^4-1}{(p^7-1)p^5}t^7 \\ & & +
\frac{p^7+p^5-p^4-1}{(p^7-1)p^4}t^6 +
\frac{p^7+p^5-p^4-1}{(p^7-1)p^3}t^5 +
\frac{p^7+p^6-p^5-1}{(p^7-1)p^3}t^4 +
\frac{p^7+p^6-p^5-1}{(p^7-1)p^2}t^3 \\ & &  +
\frac{2p^7-p^6-1}{(p^7-1)p^2}t^2 + \frac{p^8-1}{(p^7-1)p^2}t +
\frac{p^8-1}{(p^7-1)p}
\end{eqnarray*}
and
\begin{eqnarray*}
C_2(t) & = & - \frac{p-1}{(p^7-1)p} \mbox{
                }.
\end{eqnarray*}
As an illustration, we calculate the $M_i$ for $i$ in the residue
class of $3$ modulo $15$:
\begin{eqnarray}
M_{3+15e} & = & \frac{(p^7+p^6-p^5-1)p^{4+22e}}{p^7-1} -
\frac{(p-1)p^{2+15e}}{p^7-1} \label{rattermen} \\ & = & p^{4+22e}
+ \frac{(p-1)p^{9+22e}}{p^7-1} - \frac{(p-1)p^{2+15e}}{p^7-1} \nonumber \\
& = & p^{4+22e} + (p-1) \frac{p^{7e+7}-1}{p^7-1} p^{2+15e} \nonumber \\
& = & p^{4+22e} + (p-1) (p^{7e} + \cdots + p^{14} + p^7 + 1)
p^{2+15e}. \nonumber
\end{eqnarray}
Note that the two terms in (\ref{rattermen}) are not integers and
that one of them is negative. Note also that the Poincar\'e series
in the two examples are rational functions of $t$ and $p$, but
this is not the case in general.

\footnotesize{

\noindent \textsc{K.U.Leuven, Departement Wiskunde, Celestijnenlaan
200B, B-3001 Leuven, Belgium,} \\ \textsl{E-mail:}
dirk.segers@wis.kuleuven.be} \\ \textsl{URL:}
http://wis.kuleuven.be/algebra/segers/segers.htm


\begin{thebibliography}{ELM}
\bibitem[De]{Denefreport} J. Denef, \emph{Report on Igusa's local
zeta function}, S\'em. Bourbaki 741, Ast\'erisque
\textbf{201/202/203} (1991), 359-386.

\bibitem[DH]{DenefHoornaert} J. Denef and K. Hoornaert,
\emph{Newton Polyhedra and Igusa's Local Zeta Function}, J. Number
Theory \textbf{89} (2001), 31-64.

\bibitem[Hi]{Hironaka} H. Hironaka, \emph{Resolution of singularities of
an algebraic variety over a field of characteristic zero}, Ann.
Math. \textbf{79} (1964), 109-326.

\bibitem[Ig1]{Igusahigherdegree} J. Igusa, \emph{Some observations
on higher degree characters}, Amer. J. Math. \textbf{99} (1977),
393-417.

\bibitem[Ig2]{Igusabook} J. Igusa, \emph{An Introduction to the Theory
of Local Zeta Functions}, Amer. Math. Soc., Studies in Advanced
Mathematics \textbf{14}, 2000.

\bibitem[Se1]{Segersthesis} D. Segers, \emph{Smallest poles of Igusa's
and topological zeta functions and solutions of polynomial
congruences}, Ph.D. Thesis, Univ. Leuven, 2004. \\ Available on
http://wis.kuleuven.be/algebra/segers/segers.htm

\bibitem[Se2]{Segersmathz} D. Segers, \emph{On the smallest poles of
Igusa's $p$-adic zeta functions}, Math. Z. \textbf{252} (2006),
429-455.

\bibitem[Se3]{Segersmathann} D. Segers, \emph{Lower bound for the poles
of Igusa's p-adic zeta functions}, Math. Ann. \textbf{336} (2006),
659-669.

\bibitem[SV]{SegersVeys} D. Segers and W. Veys, \emph{On the
smallest poles of topological zeta functions}, Compositio Math.
\textbf{140} (2004), 130-144.
\end{thebibliography}
\end{document}